\documentclass[11pt]{article}
\usepackage[utf8]{inputenc}
\usepackage{authblk}
\usepackage{amsmath}
\usepackage{amssymb}
\usepackage{amsthm}
\usepackage{mathtools}
\usepackage{tikz}
\usepackage{listings}
\usepackage{enumitem}
\usepackage{url}
\usepackage[left=3.5cm,top=3cm,right=3cm,bottom=3cm,]{geometry}
\usepackage{pgfplots}
\usepgfplotslibrary{fillbetween}
\usepackage[
backend=biber,
style=numeric,
sorting=nyt,
maxbibnames=5
]{biblatex}

\newtheorem{thm}{Theorem}[section]
\newtheorem*{thm*}{Theorem}
\newtheorem{lemma}[thm]{Lemma}

\newtheorem{prop}[thm]{Proposition}
\newtheorem{cor}[thm]{Corollary}

\newtheorem{conj}[thm]{Conjecture}

\newtheorem*{Notation*}{Notation}

\theoremstyle{definition}
\newtheorem{defi}[thm]{Definition}
\newtheorem{ex}[thm]{Example}

\theoremstyle{remark}
\newtheorem{rem}[thm]{Remark}

\def\cross{{\mathcal C}}
\def\ehr{\mathrm{ehr}}

\DeclareMathOperator{\conv}{conv}

\DeclarePairedDelimiter\abs{\lvert}{\rvert}

\DeclarePairedDelimiter{\set}{\{}{\}}

\DeclarePairedDelimiter\floor{\lfloor}{\rfloor}

\DeclareMathOperator{\CC}{\mathbb{C}}

\DeclareMathOperator{\NN}{\mathbb{N}}

\DeclareMathOperator{\RR}{\mathbb{R}}
\DeclareMathOperator{\ZZ}{\mathbb{Z}}

\DeclareMathOperator{\triang}{\mathcal{T}}

\DeclareMathOperator{\inedge}{in}

\DeclareMathOperator{\cuts}{Cuts}
\DeclareMathOperator{\bip}{Bip}
\DeclareMathOperator{\hyptree}{ht}
\DeclareMathOperator{\CL}{CL}
\DeclareMathOperator{\lt}{lt}

\newcommand{\tildeb}[1]{\stackrel{\sim}{\smash{#1}\rule{0pt}{1.1ex}}}

\DeclareSymbolFont{extraup}{U}{zavm}{m}{n}
\DeclareMathSymbol{\vardiamondsuit}{\mathalpha}{extraup}{87}

\usepackage{xcolor}

\newcommand\restr[2]{{
  \left.\kern-\nulldelimiterspace 
  #1 
  \littletaller 
  \right|_{#2} 
  }}

\newcommand{\littletaller}{\mathchoice{\vphantom{\big|}}{}{}{}}

\title{On a Conjecture Concerning the Roots of Ehrhart Polynomials of Symmetric Edge Polytopes from Complete Multipartite Graphs}
\author{Max Kölbl\footnote{max.koelbl@ist.osaka-u.ac.jp}}
\date{}

\begin{document}

\maketitle

\begin{abstract}
    In \cite{higashitani_kummer_michalek}, Higashitani, Kummer, and Michałek pose a conjecture about the symmetric edge polytopes of complete multipartite graphs and confirm it for a number of families in the bipartite case.
    We confirm that conjecture for a number of new classes following the authors' methods and we present a more general result which suggests that the methods in their current form might not be enough to prove the conjecture in full generality.
\end{abstract}

\section*{Introduction}

A \emph{lattice polytope} is a convex polytope $P\subset\RR^n$ which can be written as the convex hull of finitely many elements of $\ZZ^n$.
Lattice polytopes arise naturally from attempts to endow combinatorial objects with a geometric structure.
A family of lattice polytopes that has garnered some attention in recent years is the \emph{symmetric edge polytope}, which is a type of graph polytope.
For graphs, we will henceforth use the notation $G=(V,E)$ where $V$ denotes the set of \emph{vertices} and $E$ denotes the set of \emph{edges} of $G$.
Given a graph $G=(V,E)$, we thus define its symmetric edge polytope as follows
\[P_G = \conv\set*{\pm(e_v-e_w)\colon \set{v,w}\in E}\subset\RR^{\abs{V}}.\]
Here, the vectors $e_v$ are elements that form a lattice basis of $\ZZ^{\abs{V}}$.
For more context on symmetric edge polytopes, see e.g. \cite{higashitani_jochemko_michalek, matsui-higashitani-nagazawa-ohsugi-hibi}.

We define the \emph{lattice-point enumerator} of a set $S\subset\RR^n$ as the function $E_S\colon\NN\to\NN$ via $E_S(k)=\abs{kS\cap\ZZ^n}$.
If $S$ is a lattice polytope, $E_S$ is a polynomial which we call the \emph{Ehrhart polynomial} of $S$.
The generating function of the Ehrhart polytope is called its \emph{Ehrhart series} and it can be written as
\[\ehr_P(t)=\sum_{k\geq 0} E_P(k) t^k = \frac{h^*(t)}{(1-t)^{d+1}},\]
where $h^*(t)$ is a polynomial with non-negative integer coefficients of degree $d$ or less.
We call this polynomial the \emph{$h^*$-polynomial} of $P$.
Both the Ehrhart polynomial and the $h^*$-polynomial hold valuable information about the underlying polytope, such as its (normalised) volume and the volume of its boundary.
A specifically remarkable piece of information encoded by the $h^*$-polynomial is that of reflexivity:
A lattice polytope is called \emph{reflexive} if its polar dual is also a lattice polytope.
By a result by Hibi \cite{hibi1992dual}, $P$ is reflexive if and only if its $h^*$-polynomial is \emph{palindromic}, i.e., $h_P^*(t)=\sum_{i=0}^d h_i^*(t)$ satisfies $h^*_i = h^*_{d-i}$ for all $0\leq i\leq d$, and its degree is equal to $\dim P$.

With some basic knowledge of generating functions (see e.g. \cite{wilf}), one can check that knowing the Ehrhart polynomial of a lattice polytope amounts to knowing its $h^*$-polynomial.
However, the converse is also true.
Given the $h^*$-polynomial $h_P^*(t)=\sum_{i=0}^d h_i^* t^i$ of some lattice polytope $P$, the Ehrhart polynomial can be written as
\[
    E_P(x) = \sum_{i=0}^d h^*_i\binom{d+x-i}{d}.
\]
For more context on Ehrhart theory, see e.g. \cite{beck-robins}.

One aspect of research in Ehrhart theory is the study of the \emph{roots} of Ehrhart polynomials when their domain and range are extended from $\NN$ to $\CC$.
For example we get \emph{Ehrhart-Macdonald reciprocity}.
\[
    E_P(-t) = (-1)^{\dim P} E_{P^\circ}(t)
\]
Here $E_{P^{\circ}}$ is the lattice-point enumerator of the \emph{interior} of $P$.
In the reflexive case, one can show from Ehrhart-Macdonald reciprocity and the palindromicity of the $h^*$-polynomial that the Ehrhart polynomial $E_P$ satisfies
\[
    (-1)^{\deg E_P} E_P(x)= E_P\left(-1-x\right).\tag{1}
\]
For the roots of such a polynomial, this means that in addition to the symmetry across the real axis, there is a second symmetry across the \emph{canonical line}, i.e., the set
\[
    \CL=\set*{z\in\CC\colon\Re(z)=-\frac{1}{2}}
\]
where $\Re(z)$ denotes the real part of $z\in\CC$.
Thus, it is natural to ask, what kind of polytopes have all of their Ehrhart polynomial roots \emph{on} $\CL$.
First steps in this direction were made in \cite{bump-choi-kurlberg-vaaler} and \cite{rodriguez-villegas}, albeit in different contexts.
In \cite{matsui-higashitani-nagazawa-ohsugi-hibi}, the study of \emph{$\CL$-polytopes}, i.e., polytopes with all their Ehrhart polynomial roots on $\CL$, has been initiated as a field of study in its own right.
For low dimensions, a full classification was found in \cite{hegedus-higashitani-kasprzyk}.
Some known classes include cross-polytopes (which we will define properly later), standard reflexive simplices, and root polytopes of type A.

In \cite{higashitani_kummer_michalek}, the authors studied the roots of the Ehrhart polynomials of symmetric edge polytopes of the complete bipartite graphs $K_{2,n}$ and $K_{3,n}$ and were able to prove that both these classes are $\CL$-polytopes.
This extends the case of cross-polytopes, which are the symmetric edge polytopes of $K_{1,n}$.
They accomplished that by using the technique of \emph{interlacing polynomials}, i.e., polynomials whose roots alternate on a given totally ordered set.
For an in-depth treatment of the theory of interlacing polynomials, see \cite{fisk2006polynomials}.
The authors gave the following conjecture.

\begin{conj}[Conjecture 4.1 in \cite{higashitani_kummer_michalek}]\label{conjecture: HKM}
    \begin{itemize}
        \item[(i)] For any complete multipartite graph $K_{a_1,\dots,a_k}$ the Ehrhart polynomial $E_{a_1,\dots,a_k}$ has its roots on $\CL$.
        \item[(ii)] Suppose $a_1\leq\dots\leq a_k$. Any two Ehrhart polynomials $E_{a_1,\dots,a_k}$ and $E_{a_1,a_2,\dots,a_k-1}$ interlace on $\CL$.
    \end{itemize}
\end{conj}

We are able to prove $\CL$-ness of the symmetric edge polytopes of $K_{1,1,n},K_{1,2,n}, K_{1,1,1,n}$, as well as some conditional results.

The structure of the paper is as follows.
In Section 2, we will recall a number of known results about symmetric edge polytopes and interlacing polynomials.
In Section 3, we will compute the $h^*$-polynomials of symmetric edge polytopes of complete tripartite graphs (Theorem \ref{theorem: tripartite h*}).
We will do that by following the methods used by Higashitani, Jochemko, and Michałek in \cite{higashitani_jochemko_michalek} who computed the $h^*$-polynomials of symmetric edge polytopes of complete bipartite graphs by finding a half-open triangulation using Gröbner basis techniques.
In particular, in Theorem \ref{theorem: gröbner basis}, we find a cubic Gröbner basis for every symmetric edge polytope of multipartite graphs, which is a reduction of the Gröbner basis from \cite{higashitani_jochemko_michalek}.
The reduced Gröbner basis assumes a certain ordering of the edges of the graph, which the unreduced one does not.
However, we will also show that no quadratic reduction is possible if the underlying graph contains $K_{2,2,2}$ as a subgraph.

Section 4 contains the new recursive relations in the style of the ones found by Higashitani, Kummer, and Michałek in \cite{higashitani_kummer_michalek} and in Theorem \ref{theorem: new relations}, we provide new evidence for Conjecture \ref{conjecture: HKM}.
In Section 5 we will investigate a connection between the $\gamma$-vector of the $h^*$-polynomial of an Ehrhart polynomial and the existence of recursive relations that can be used to prove interlacing.
In particular, in Theorem \ref{theorem: recursive relations general} and Corollaries \ref{corollary: recursion cor 1} and \ref{corollary: recursion cor 2}, we show that the type of recursive relations in \cite{higashitani_kummer_michalek} and Proposition \ref{proposition: new relations} can be found for arbitrary complete bipartite graphs.

\section*{Acknowledgements}
I would like to express my gratitude to my advisor Akihiro Higashitani and to Rodica Dinu for many useful discussions and invaluable advice.

\section{Preliminaries}

\subsection{Symmetric edge polytopes and triangulations}

The content of this subsection is based on \cite{higashitani_jochemko_michalek}.
Our goal is to build the machinery necessary to compute $h^*$-polynomials of symmetric edge polytopes.
The first step in this direction is to understand their facets.

\begin{prop}[Theorem 3.1 in \cite{higashitani_jochemko_michalek}]
    Let $G = (V, E)$ be a finite simple connected graph.
    Then $\lambda\colon V\to\ZZ$ is facet defining if and only if
    \begin{itemize}
        \item[(i)] for any edge $e=\set{u,v}$ we have $\abs{\lambda(u)-\lambda(v)}\leq 1$, and
        \item[(ii)]  the subset of edges $E_\lambda=\set{\set{u,v}\in E\colon\abs{\lambda(u)-\lambda(v)}=1}$ forms a spanning subgraph of $G$.
    \end{itemize}
\end{prop}

With this, we can restrict ourselves to the case of complete multipartite graphs.
The bipartite case and the $\geq 3$-partite case work slightly differently, but we will only cite the latter case.
For the bipartite case, see Proposition 3.4 in \cite{higashitani_jochemko_michalek}.

\begin{prop}[Proposition 3.5 in \cite{higashitani_jochemko_michalek}]\label{proposition: HJM facets}
     Let $k\geq 3$ and $G = K_{a_1,\dots,a_k}$ be a complete $k$-partite graph with vertex set $V = \bigsqcup_{i=1}^k A_i$ and edges between any two vertices if and only if they belong to different $A_i$.
     Then $\lambda\colon V\to\ZZ$ is facet defining if and only if $\lambda$, up to a constant, satisfies one of the following conditions.
     \begin{itemize}
        \item[(i)] $\lambda(A_i)=\set{-1,1}$ for some $1\leq i\leq k$ and $\restr{\lambda}{A_j}=0$ for all $i\neq j$, or
        \item[(ii)] $\lambda(V)=\set{0,1}$ and 
        \begin{itemize}
            \item[(a)] $\restr{\lambda}{A_i}$ is constant for every $A_i$, or
            \item[(b)] there exists an $i$ such that $\lambda(A_i)=\set{0,1}=\lambda(\bigcup_{j=1}^k A_j\setminus A_i)$.
        \end{itemize}
    \end{itemize}
    In particular, the polytope symmetric edge polytope of $G$ has $2^{\sum^k_{i=1} a_i}-\sum^k_{i=1}(2a_i-2)-2$ facets.
\end{prop}

The next step is to find a \emph{triangulation} of $P_G$.
A triangulation of a $d$-dimensional polytope $P$ is a decomposition of $P$ into simplices such that the intersection of every pair of simplices is a simplex of lower dimension.
A simplex is called \emph{unimodular} if its vertices span $\ZZ^d$.
A triangulation consisting of unimodular simplices is called a \emph{unimodular triangulation}.
In order to get that, we employ methods from commutative algebra.
The following result details a \emph{Gröbner basis} of $P_G$ if $G$ is a complete multipartite graph.
We start by picking an orientation of $G$ and for every edge $e$ in $G$, we consider the variables $x_e$ and $y_e$ which represent directed versions of $e$ where $x_e$ follows the orientation and $y_e$ goes against it.
If there is a second, independent orientation on some of the edges, we will write $p_e$ for the variable representing the directed edge following the \emph{second} orientation and $q_e$ for the opposite edge.
In particular, for every edge, we have $\set{x_e,y_e}=\set{p_e,y_e}$.

\begin{prop}[Proposition 3.8 in \cite{higashitani_jochemko_michalek}]\label{proposition: HJM GB}
    Let $z < x_{e_1} < y_{e_1} < \dots < x_{e_k} < y_{e_k}$ be an order on the edges.
    Then the following collection of three types of binomials forms a Gröbner basis of the toric ideal of the symmetric edge polytope of $G$ with respect to the degrevlex order:
    \begin{itemize}
        \item[(1)] For every $2k$-cycle $C$, with fixed orientation, and any $k$-element subset $I$ of edges of $C$ not containing the smallest edge
        \[\prod_{e\in I}p_e -\prod_{e\in C\setminus I} q_e.\]
        \item[(2)] For every $(2k + 1)$-cycle $C$, with fixed orientation, and any $(k + 1)$-element subset $I$ of edges of $C$
        \[\prod_{e\in I} p_e - z\,\prod_{e\in C\setminus I} q_e.\]
        \item[(3)] For any edge $e$
        \[x_e y_e - z^2.\]
    \end{itemize}
    The leading monomial is always chosen to have positive sign.
\end{prop}

This Gröbner basis has a geometric interpretation.
Given a symmetric edge polytope $P_G$ for some graph $G=(V,E)$, if $\Delta$ is a simplex of a unimodular triangulation of $P_G$, a vertex of $\Delta$ is either a vertex of $P_G$ or its unique interior lattice-point.
Since we can identify every vertex of $P_G$ with a directed edge of $G$, we can further identify it with the corresponding variable $x_e$ or $y_e$.
The remaining variable, $z$, can be identified with the interior lattice-point.
In this way, every square-free degree $d+1$ monomial in these variables corresponds to the vertices of a $d$-dimensional lattice simplex $\Delta\subset P_G$, provided the vertices are affinely independent.
In this setup, let $T$ be the set of square-free degree $d$ monomials $m$ in the variables $\set{x_e,y_e\colon e\in E}$ such that there exists no Gröbner basis element $p$ whose leading term $\lt(p)$ divides $m$. 
Then $T$ corresponds to a unimodular triangulation of the boundary of $P_G$.
In particular, for every $m$ the set of corresponding edges of $G$ induces a directed spanning tree of $G$.
Thus we can define the set $\triang$ which consists of all spanning trees of $G$ that come from monomials in $T$.
Now we fix a vertex $r$ in $G$.
Given a directed spanning tree $S\in\triang$, we call an edge $e$ in $S$ \emph{ingoing} if the unique path starting at the foot of $e$ and ending in $r$, the path includes $e$.
Otherwise we call it \emph{outgoing}.
We denote the number of ingoing edges of $S$ by $\inedge(S)$.
With all this information, one can compute the $h^*$-polynomial of $P_G$.

\begin{prop}[Proposition 4.6 in \cite{higashitani_jochemko_michalek}]
    Let $h^*_G(t)=\sum_{i=0}^d h^*_i t^i$.
    Then
    \[ h_i^* = \abs{\set{T\in\triang\colon\inedge (T) = i}}. \]
\end{prop}

The authors use this formula to compute the $h^*$-polynomial in the case of complete bipartite graphs.

\begin{prop}[Theorem 4.1 in \cite{higashitani_jochemko_michalek}]\label{proposition: HJM h*}
    For all $a, b \geq 0$ let $h^*_{a,b}(t)$ denote the $h^*$-polynomial of the symmetric edge polytope of $K_{a+1,b+1}$.
    Then
    \[
        h_{a,b}^*(t) = \sum_{i=0}^{\min\set{a,b}}\binom{2i}{i}\binom{a}{i}\binom{b}{i}t^i(1+t)^{a+b+1-2i}.
    \]
\end{prop}

\subsection{Symmetric edge polytopes and graph constructions}

The content of this subsection is based on \cite{ohsugi_tsuchiya}.
We will revisit two useful results pertaining to the computation of $h^*$-polynomials under certain graph constructions.

The first one out of the two needs some preparation.
A \emph{hypergraph} $\mathcal{H}$ is a set $V$ and a set $E$ of nonempty subsets of $V$ called \emph{hyperedges}.
We can associate a bipartite graph $\bip(\mathcal{H})$ to $\mathcal{H}$ whose bipartite classes are given by elements of $V$ and $E$ respectively with an edge between a $v\in V$ and an $e\in E$ if $v\in e$.
A \emph{hypertree} is a function $f\colon E\to \set{0,1,\dots}$ such that there exists a spanning tree $\Gamma$ of $\bip(\mathcal{H})$ whose vertices $e\in E$ have degree $f(e)+1$.
In this case, we say that $\Gamma$ \emph{induces} $f$.
The set of hypertrees of $\mathcal{H}$ shall be denoted by $\hyptree(\mathcal{H})$.
Let us now fix a total order of $E$.
A hyperedge $e$ is called \emph{internally active} with respect to $f$ is there exists no $e^\prime<e$ such that increasing $f(e^\prime)$ by $1$ and decreasing $f(e)$ results in a different hypertree.
A hyperedge that is not internally active with respect to $f$ is called \emph{internally inactive}.
We denote the number of internally inactive edges of $f$ by $\iota(f)$.
The \emph{interior polynomial} $I_{\mathcal{H}}$ of $\mathcal{H}$ is then defined as
\[
    I_{\mathcal{H}}(t)=\sum_{f\in\hyptree(\mathcal{H})}t^{\iota(f)}.
\]
Given a graph $G=\bip(\mathcal{H})$ for some hypergraph $\mathcal{H}$, we define $I_G = I_{\mathcal{H}}$.

Next, recall that a \emph{cut} of a graph $G=(V,E)$ is a subgraph $G_C=(V,E_C)$ of $G$ where $C$ is a subset of $V$ and $E_C\subset E$ is the set of edges of $G$ with one end in $C$ and one end not in $C$.
We denote the set of cuts of $G$ by $\cuts(G)$.

Lastly, we need two special graph constructions.
Given a graph $G$ with vertex set $[d]$, let $\widehat{G}$ describe the the \emph{suspension} of $G$,
i.e., the graph on the set $[d+1]$ with the same edge set as $G$ but with the vertex $d+1$ connected to all the others.
If $G$ is bipartite with bipartite classes $V$ and $W$, its \emph{joint bipartite suspension} $\tildeb{G}$ is the graph on $[d+2]$ such that $d+1$ connects to all the edges in $V$, $d+2$ connects to all the edges in $W$ and $d+1$ and $d+2$ connect to each other.
Like this, we can cite the following theorem.

\begin{prop}[Theorem 4.3 in \cite{ohsugi_tsuchiya}]\label{proposition: OT suspension}
    Let $G$ be a finite graph on the vertex set $[d]$.
    Then the symmetric edge polytope of $\widehat{G}$ is unimodularly equivalent to a reflexive polytope whose $h^*$-polynomial is
    \[
        h^*_{\widehat{G}}(t) = (1+t)^df_G\left(\frac{4t}{(1+t)^2}\right),
    \]
    where $f_G(t) = \frac{1}{2^{d-1}}\sum_{H\in\cuts(G)} I_{\tildeb{H}}(t)$.
\end{prop}

We need this Proposition for a specific family of symmetric edge polytopes, namely the one that comes from complete multipartite graphs of the form $K_{1,1,1,n}$, which we shall treat in the following example.

\begin{ex}\label{example: h*_1,1,1,n}
    The graph $G=K_{1,1,n}$ gives rise via suspension to the graph $\widehat{G}=K_{1,1,1,n}$.
    We denote the the vertices in the first two tripartite classes of $G$ by $a$ and $b$ respectively.
    The remaining vertices shall be denoted by the integers $1,\dots,n$.
    First, we need to understand the cuts $G_C$ of $G$.
    There are two primary types: one type where without loss of generality $a,b\not\in C$, and one type where $a\in C$, $b\not\in C$.
    Thus, for every subset $S\subseteq [n]$ we get a cut set $C_1=S$ and a cut set $C_2=S\cup\set{a}$.
    Assume $\abs{S}=m$.
    Now we need to understand the hypertrees associated to the joint bipartite suspensions $\tildeb{G}_{C_i}$ for $i=1,2$, which we regard as functions $f\colon \set{n+3}\cup C_i\to\set{0,1,\dots}$.
    In the order of the hyperedges, it is convenient to regard $n+3$ as the smallest edge.
    This way, it can never be an internally inactive edge and we can focus on the elements of $C_i$ instead.
    One can check that a hyperedge $e\in C_i$ is internally inactive if and only if $f(e)>0$.
    Without loss of generality, we can assume that every inducing spanning tree of a hypertree contains the edge $\set{n+3,n+4}$ and for every $c\in C_i$ the edge $\set{n+4,c}$.
    From here, one can check that the interior polynomial of $\tildeb{G}_{C_1}$ is $\binom{m}{2}t^2+2mt+1$ and that of $\tildeb{G}_{C_2}$ is $m(n-m)t^2+(n+1)t+1$.
    Summing up, we get
    \begin{align*}
        f_G(t) &= \frac{1}{2^{n+1}}\sum_{i=0}^n\binom{n}{m}\left(\left(\binom{m}{2}+m(n-m)\right)t^2 + (2m+n+1)t + 2\right) \\
        &= \frac{3(n-1)n}{2^4}t^2 + \frac{2n+1}{2}t + 1.
    \end{align*}
    Thus, we obtain the $h^*$-polynomial of the symmetric edge polytope of $K_{1,1,1,n}$:
    \[
        h^*_{1,1,1,n}(t) = 3(n-1)n(1+t)^{n-2}t^2 + 2(2n+1)(1+t)^{n}t + (1+t)^{n+2}
    \]
\end{ex}

The other very useful theorem is the following.

\begin{prop}[Proposition 4.4 in \cite{ohsugi_tsuchiya}]\label{proposition: OT contraction}
    Let $G$ be a bipartite graph on $[d]$ and let $e$ be an edge of $G$.
    Then we have
    \[
        h_{G}^*(t)=(t+1)h_{G/e}^*(t)
    \]
    where $G/e$ denotes the $G$ after contraction of the edge $e$.
\end{prop}

\subsection{Interlacing}

In this subsection, we will recall some basic results concerning the theory of interlacing polynomials.
We start by defining what exactly we mean by ``interlacing''.

\begin{defi}
    Let $(R,\preceq)$ be a totally ordered subset of $\CC$ and let $f,g$ be polynomials of degree $d+1$ and $d$ respectively whose roots lie in $R$.
    Then $g$ \emph{$R$-interlaces} $f$ if
    \[
    a_1\preceq b_1\preceq a_2\preceq b_2\preceq \dots\preceq b_d\preceq a_{d+1}
    \]
    where the $a_i$ are the roots of $f$ and the $b_i$ are the roots of $g$.
\end{defi}

Going forward, we will only consider Ehrhart polynomials of polytopes whose roots lie on $\CL$ and we will focus on $\CL$-interlacing.
Hence, we give the following proposition in a much more specialised form than the one found in \cite{higashitani_kummer_michalek}.

\begin{prop}[Lemmas 2.3, 2.4, 2.5 in \cite{higashitani_kummer_michalek}]\label{proposition: HKM interlacing}
    Let $f,g,h_1,\dots,h_n$ be real polynomial such that $\deg f =\deg g+1 = \deg h_i+2$ for all $1\leq i\leq n$ which all satisfy the symmetry relation in Equation (1) as given in the introduction.
    Assume the identity
    \[
    f(x) = (2x+1)\alpha g(x) + \sum_{i=1}^n \alpha_i h_i(x)
    \]
    where $\alpha,\alpha_i>0$ for all $i$.
    Then $\sum_{i=1}^n \alpha_i h_i$ $\CL$-interlaces $g$ if for every $i$, $h_i$ $\CL$-interlaces $g$.
    Also, the following are equivalent.
    \begin{itemize}
        \item[(a)] $\sum_{i=1}^n \alpha_i h_i$ $\CL$-interlaces $g$,
        \item[(b)] $g$ $\CL$-interlaces $f$.
    \end{itemize}
    If this is the case, $(2x+1)\sum_{i=1}^n \alpha_i h_i$ $\CL$-interlaces $f$.
\end{prop}

An important class of reflexive polytopes is the class of \emph{cross-polytopes} which are defined as the convex hull of the vectors $\pm e_1,\pm e_2,\dots,\pm e_n\in\RR^n$.
In particular, they are unimodularly equivalent to symmetric edge polytopes of trees with $n$ edges.
The Ehrhart polynomial of the $n$-dimensional cross-polytope (the $n$-th \emph{cross-polynomial}) is given by
\[
\cross_{n}(x) = \sum_{k=0}^n \binom{n}{k}\binom{n+x-k}{n}.
\]
Cross-polynomials are the first class of examples to showcase the usefulness of Proposition \ref{proposition: HKM interlacing}.

\begin{prop}[Example 3.3 in \cite{higashitani_kummer_michalek}]\label{proposition: cross-polynomial relations}
    For any $n\geq 2$, cross-polynomials satisfy the recursive relation
    \[
        \cross_n(x) = \frac{1}{n}(2x+1)\cross_{n-1}(x) + \frac{n-1}{n}\cross_{n-2}(x).
    \]
\end{prop}

Other classes of examples were found by Higashitani, Kummer, and Michałek in \cite{higashitani_kummer_michalek}.
The authors found three recursive relations among Ehrhart polynomials $E_{1,n},E_{2,n},E_{3,n}$ of the symmetric edge polytopes from the complete bipartite graphs $K_{1,n}, K_{2,n}, K_{3,n}$.

\begin{prop}[Proposition 4.5 in \cite{higashitani_kummer_michalek}]\label{proposition: HKM known relations}
    The following relations hold:
    \begin{align*}
        E_{2,n}(x) &= \frac{1}{2} (2x+1) E_{1,n}(x) + \frac{1}{2} E_{1,n-1}(x), \\
        E_{2,n}(x) &= \frac{1}{n} (2x+1) E_{2,n-1}(x) + \frac{1}{2n}\left(n E_{1,n-1} + (n-2)(2x+1) E_{1,n-2}\right), \\
        E_{3,n+1}(x) &= \frac{(2x+1)(3n^2 + 13n + 16)}{8(n^2 + 5n + 6)} E_{2,n+1} \\
        &+ \frac{n^3 13n^2 + 18n}{8(n - 1) (n^2 + 5n + 6)} E_{2,n} + \frac{4n^3 + 9n^2 - 13n - 32}{8(n - 1) (n^2 + 5n + 6)} E_{1,n+1}.
    \end{align*}
\end{prop}

This leads to the following proposition.

\begin{prop}[Lemmas 4.6-4.8, Theorem 4.9 in \cite{higashitani_kummer_michalek}]\label{proposition: HKM known interlacings}
    The following statements hold.
    In the following, $E_{a,b}$ denotes the Ehrhart polynomial of $P_{K_{a,b}}$.
    \begin{itemize}
        \item[(a)] For every $n\geq 1$, $E_{1,n}$ $\CL$-interlaces $E_{1,n+1}$.
        \item[(b)] For every $n\geq 1$, the Ehrhart polynomials of $K_{1,n}$ and $(2k+1)K_{1,n-1}$ $\CL$-interlace $E_{2,n}$.
        \item[(c)] For every $n\geq 1$, $E_{2,n}$ $\CL$-interlaces $E_{2,n+1}$.
        \item[(d)] For every $n\geq 1$, $E_{2,n}$ $\CL$-interlaces $E_{3,n}$.
    \end{itemize}
    In particular, for every $n\geq 1$ the Ehrhart polynomials of $K_{m,n}$ is a $\CL$-polynomial if $m\leq 2$.
\end{prop}

\section{A Reduced Gröbner Basis}

We start by describing an edge order.
First we denote the multipartite classes of vertices of $K_{a_1,\dots,a_k}$ by $A_1, A_2,\dots,A_k$ and then we pick an order of the vertices which satisfies the following condition.
If $v\in A_i$ and $w\in A_j$, then $v<w$ if and only if $i<j$.
Let $e=\set{v,w}$ and $e^\prime=\set{v^\prime,w^\prime}$ be edges in $K_{a_1,\dots,a_k}$.
Without loss of generality, we may assume $v<w$ and $v^\prime < w^\prime$.
Then $e<e^\prime$ if and only if $v < v^\prime$ or $v=v^\prime$ and $w<w^\prime$.

Let $a,b$ be vertices with an edge between them.
We will denote by $x_{a,b}$ the directed edge from $a$ to $b$ and by $x_{b,a}$ the edge going the other way.
The variable which corresponds to the unique interior lattice-point of $P_{K_{a_1,\dots,a_k}}$ will still be denoted by $z$.

\begin{thm}\label{theorem: gröbner basis}
    With the described edge order, the Gröbner basis from Proposition \ref{proposition: HJM GB} is at most cubic for every complete multipartite graph $K_{a_1,a_2,\dots,a_k}$.
    The elements of the reduced Gröbner basis are of the following form.
    \begin{itemize}
        \item [(1)] Let $a\in A_i$ and $b\in A_j$ with $i\neq j$.
        Then the following polynomial is a Gröbner basis element.
        \[ x_{a,b} x_{b,a} - z^2 \]
        \item[(2)] Let $a\in A_i, b\in A_j, c\in A_\ell$ with $i,j,\ell$ all different.
        Then the following polynomial is a Gröbner basis element.
        \[ x_{a,b} x_{b,c} - z x_{a,c} \]
        \item[(3)] Let $a,b,c,d$ be vertices such that the edges $\set{a,b},\set{b,c},\set{c,d},\set{a,d}$ all exist and $a$ is the smallest vertex.
        Then the following polynomial is a Gröbner basis element if and only if $b$ and $d$ lie in the same $A_i$.
        \[ x_{b,c} x_{c,d} - x_{b,a} x_{a,d} \]
        We call these polynomials \emph{Gröbner basis elements of type (3a).}
        Furthermore, the following polynomial is a Gröbner basis element if and only if $b<d$.
        \[ x_{b,c} x_{d,a} - x_{b,a} x_{d,c} \]
        We call these polynomials \emph{Gröbner basis elements of type (3b).}
        In particular, $a,b,c,d$ lie across either $2$, $3$, or $4$ multipartite classes.
        \item[(4)] Let $a,b,c,d,e$ be vertices such that the edges $\set{a,b},\set{b,c},\set{c,d},\set{d,e},\set{a,e}$ all exist.
        Then the following polynomial is a Gröbner basis element if and only if $a,b,c\in A_1\cup A_2$, $a,c$ lie in the same $A_i$, and $b$ is the smallest vertex in $A_1$ or $A_2$.
        \[ x_{a,b} x_{b,c} x_{d,e} - z\, x_{d,c} x_{a,e} \]
        In particular, $a,b,c,d,e$ lie across either $3$ or $4$ multipartite classes.
        \item[(5)] Let $a,b,c,d,e,f$ be vertices such that the edges $\set{a,b}$,$\set{b,c}$,$\set{c,d}$,$\set{d,e}$,$\set{e,f}$,$\set{a,f}$ all exist.
        Then the following polynomial is a Gröbner basis element if and only if
        \begin{itemize}
            \item[(i)] $c$ and $f$ lie in the same $A_i$,
            \item[(ii)] $b>d$, or $b$ and $e$ lie in the same $A_i$, or $c<e$,
            \item[(iii)] $a$ and $d$ lie in the same $A_i$, or $f<d$.
        \end{itemize}
        \[ x_{a,b} x_{b,c} x_{d,e} - z\, x_{d,c} x_{a,e} \]
        In particular, $a,b,c,d,e$ lie across either $3$, $4$, or $5$ multipartite classes.
    \end{itemize}

    More generally, for every complete multipartite graph $K_{a_1,\dots,a_k}$ which contains $K_{2,2,2}$ as a subgraph, the Gröbner basis in Proposition \ref{proposition: HJM GB} has an element of degree $3$ regardless of the edge order.
\end{thm}

\begin{proof}
    One can check that all the listed elements indeed come from directed cycles in the way described in Proposition \ref{proposition: HJM GB}.
    To check the reducedness of a Gröbner basis element $p$, it is enough to find another element $q$ of lower degree such that $\lt (q) | \lt (p)$, where $\lt(p)$ and $\lt(q)$ are the leading terms of $p$ and $q$ respectively.
    Since all elements are of degree at least $2$, we can see that (1), (2), and (3) are indeed  not redundant.
    For (4) and (5) we may notice that the given restrictions correspond to indivisibility of the polynomials by the leading terms of Gröbner basis elements of type (2) or (3).

    Now we can go on to show that no further elements are contained in the Gröbner basis.
    Firstly, let $C$ be a directed cycle of length $7$ or greater.
    Assume the set $I$ which defines the leading term of the polynomial $p_{C,I}$ contains two adjacent directed edges $(a,b)$ and $(b,c)$.
    One can check that there exists an element of type (2) or (3) whose leading term contains these edges unless $a$ and $c$ both lie in $A_i$ and $b$ is the smallest vertex in $A_j$ with $\set{i,j}=\set{1,2}$.
    In this case, $C$ cannot be even because $I$ contains the smallest edge.
    That means that there exists a vertex $d$ not in $A_1$ or $A_2$.
    Further, we assume that $(a,b)$ and $(b,c)$ are the only pair of adjacent edges in $C$.
    Due to the size constraint of $I$, every vertex other than $b$ is part of one directed edge in $I$.
    Thus, there exists a directed edge $(d,e)$ or $(e,d)$ whose associated variable forms the leading term of a Gröbner basis element of type (4) together with $x_{a,b}$ and $x_{b,c}$.

    Next, let $C$ be a directed cycle of length $8$ or greater.
    We assume that the set $I$ which defines the leading term of the polynomial $p_{C,I}$ contains no adjacent directed edge.
    Thus, $C$ is necessarily an even cycle.
    We denote the vertices of the cycle in order by $a_0,b_0,a_1,b_1,\dots,a_n,b_n$ such that $a_0$ is the smallest vertex and we can assume up to orientation that $(b_i,a_{i+1})\in I$ for $0\leq i\leq n-1$ and $(b_n,a_0)\in I$.
    If $a_1$ and $b_i$ for $i>1$ lie in the same multipartite class, we get a smaller cycle $C^\prime$ containing $a_0, b_0, a_1, b_i,a_{i+1},\dots,b_n$ with $I^\prime\subset I$ such that $\lt (p_{C^\prime,I^\prime}) | \lt(p_{C,I})$.
    Thus, we may assume that $a_1$ and $b_i$ for $i>1$ all lie in the same multipartite class.
    This puts $a_2$ in a different class from $b_n$.
    As a consequence, the directed cycle $C^\prime$ on the vertices $a_0,b_0,a_1,b_1,a_2,b_n$ with $I^\prime=\set{(b_0,a_1),(b_1,a_2),(b_n,a_0)}$ yields a polynomial $p_{C^\prime,I^\prime}$ whose leading term divides that of $p_{C,I}$.
    
    Lastly, we prove the second part of the theorem.
    Let $K_{2,2,2}$ be a subgraph of $K_{a_1,\dots,a_k}$.
    Then there exists a directed $6$-cycle with vertices $a,b,c,d,e,f$ with edges 
    \[(a,b),(b,c),(c,d),(d,e),(e,f),(f,a)\]
    such that $(a,b)$ is the smallest edge of $K_{2,2,2}$.
    This gives rise to the polynomial $c = x_{b,c} x_{d,e} x_{f,a} - x_{b,a} x_{d,c} x_{f,e} $ which is an element of the Gröbner basis.
    We can verify that there does not exist another Gröbner basis element whose leading monomial divides that of $c$.
\end{proof}

If we want to use this Gröbner basis to find a unimodular triangulation, we may notice that not all elements need to be considered.
We know that for every unimodular simplex in the triangulation, its vertices that lie in the boundary of $P_G$ all lie within the same facet.
Further we know by Proposition \ref{proposition: HJM facets} that these facets are given by labelings of the vertices of $K_{a_1,\dots,a_k}$ which satisfy specific conditions.
Indeed, edge configurations induced by Gröbner basis elements of types (1), (2), and (4) do not occur in any facet-inducing spanning subgraph.
Configurations induced by elements of type (5) only appear in spanning trees of type (ii).
Among the configurations induced by elements of type (3), both varieties appear in facet-inducing spanning subgraphs of type (i), whereas in type (ii) only type (3b) elements appear.

\subsection{Simplices in type (i) facets}

To start, we will establish some terminology and notation.
Let $G=(V,E)$ be a graph and $P_G$ its symmetric edge polytope.
As seen in Proposition \ref{proposition: HJM facets}, a facet of $P_G$ is induced by an integer valued function on $V$.
We will henceforth call such a function a \emph{labeling} on $V$.
We will denote it with a lowercase Greek letter such as $\lambda$.
Following this, we call a vertex $v$ \emph{$\ell$-labled} if $\lambda(v)=\ell$.
The facet of $P_G$ induced by $\lambda$ shall be denoted by $\mathcal{F}_\lambda$.
The spanning subgraph of $G$ induced by $\lambda$ shall be denoted by $\restr{G}{\lambda}$.
The simplices in the unimodular triangulation of $P_G$ will be denoted by the symbol $\Delta$ and the associated directed spanning tree by $T_\Delta$.
For the unimodular triangulation itself, we will write $\triang$.
Given a labeling $\lambda$, the set $\triang_\lambda\subset\lambda$ is the set of simplices which lie in $\mathcal{F}_\lambda$.
Lastly, we define the set $\triang_{(i)}$ as the union of all the $\triang_\lambda$ where $\lambda$ is a type (i) facet, and the set $\triang_{(ii)}$ analogously.

The following definition should be viewed with an eye toward Gröbner basis elements of type (3b):
Let $\lambda$ be a facet-inducing labeling and let $A$ and $B$ be sets of vertices such that no $a\in A$ lies in the same multipartite class as a $b\in B$.
Further, assume that $\lambda(a)=\lambda(b)-1$ for every $a\in A$ and $b\in B$, and that $\restr{\lambda}{A}$ and $\restr{\lambda}{B}$ are constant.
The spanning subgraph corresponding to this situation would contain a directed edge from every element of $A$ to every element of $B$.
The following definition tells us which subsets of edges from $A$ to $B$ can be included in ``legal'' spanning trees with respect to the Gröbner basis from Theorem \ref{theorem: gröbner basis}. 

\begin{defi}
    Let $A$ and $B$ be disjoint finite totally ordered sets.
    A \emph{planar spanning tree} between $A$ and $B$ is a subset $E$ of $A\times B$ such that
    \begin{itemize}
        \item[(i)] $\abs{E}=\abs{A\cup B}-1$,
        \item[(ii)] every element of $A$ and $B$ is contained in at least one element of $E$,
        \item[(iii)] if $(a,b)$ and $(a^\prime,b^\prime)$ are elements of $E$, then $a<a^\prime$ implies $b<b^\prime$. 
    \end{itemize}
    The number of planar spanning trees is $\binom{a+b-2}{b-1}$ where $a$ and $b$ are the cardinalities of $A$ and $B$ respectively.
\end{defi}

\begin{prop}
    Let $K_{a_1,a_2,\ldots,a_k}$ be a complete multipartite graph with multipartite classes of vertices $A_1,A_2,\dots,A_k$ and let $P_{K_{a_1,a_2,\ldots,a_k}}$ be its associated symmetric edge polytope.
    
    Then the polynomial $h^{(i)}_{a_1,a_2,\ldots,a_k}=\sum_{\Delta\in\triang_{(i)}} t^{\inedge(T_\Delta)}$ is given by
    \begin{align*}
        h^{(i)}_{a_1,a_2,\ldots,a_k}(t) &= \sum_{i=0}^{a-a_i-1} \sum_{j=1}^{a_1-1} p(a,a_1,i,j)\binom{a_1+j-i-2}{j-1}\left(t^{i+j+1}+t^{a-i-j-2}\right) \\
        +\sum_{m=2}^k &\sum_{i=0}^{a-a_m-1} \sum_{j=1}^{a_m-1} p(a,a_m,i,j)\binom{a-a_m+j-i-2}{a-a_m-i-1}\left(t^{i+j}+t^{a-i-j-1}\right)
    \end{align*}
    where $\inedge(T)$ is the number of ingoing edges of $T$, $a=a_1+a_2+\dots+a_k$ and
    \[p(x,y,i,j) = \binom{x-y-1}{i}\binom{y-1}{j}\binom{y+i-j-1}{i}.\]
\end{prop}

\begin{proof}
    We fix a labeling $\lambda\colon\bigsqcup_{j=1}^k A_j\to\set{-1,0,1}$ corresponding to a facet of type (i).
    This means that for one $A_m$ we get $\restr{\lambda}{A_m}=\set{-1,1}$ and all the remaining vertices are mapped to $0$.
    For the graph $\restr{K_{a_1,a_2,\ldots,a_k}}{\lambda}$ this means that every vertex $v\not\in A_m$ has an edge $(v,w)$ if $w$ is $1$-labeled and an edge $(w,v)$ if $w$ is $-1$-labeled.
    On the other hand, all the vertices in $A_m$ only have edges leading into them or out of them, depending on their labeling.
    Consider now a spanning tree $T_\Delta$ with $\Delta\in\triang_\lambda$.
    Within $T_\Delta$, the Gröbner basis elements of type (3a) block every vertex $0$-labeled vertex $v$ (with one exception) from having edges of the form $(v,w)$ and $(w,v)$ at the same time.
    The exception is the smallest $0$-labeled vertex in the graph, which we will denote by $v_0$.
    Thus, we obtain two subsets of $\bigsqcup A_j\setminus A_m$:
    the subset $P$ of vertices $v$ whose edges are of the form $(v,w)$, and
    the subset $N$ of vertices $v$ whose edges are of the form $(w,v)$.
    In particular, $P\cap N =\set{v_0}$.
    With these conditions $A_m$ naturally splits into two disjoint subsets $A_m^+=\set{v\in A_m\colon \lambda(v)=1}$ and $A_m^-=\set{v\in A_m\colon \lambda(v)=-1}$.
    Taking the Gröbner basis elements of type (3b) into account, $T_\Delta$ is the disjoint union of a planar spanning tree between $P$ and $A_m^+$ and a planar spanning tree between $N$ and $A_m^-$.

    Next, we want to count the number of ingoing edges.
    Let $r$ denote the smallest vertex in $\bigsqcup A_j$ and let $v$ be some element in $P$ different from $r$.
    The edge $e$ containing $v$ in the unique path from $r$ to $v$ is ingoing.
    In a similar way, if $v$ is any element in $A_m^-$ different from $r$, the edge containing $v$ in the unique path from $r$ to $v$ is also ingoing.
    Every other edge is outgoing.
    We get a total of four cases:
    \begin{itemize}
        \item[(a)] $m=1$ and $\lambda(\min A_m)=1$,
        \item[(b)] $m=1$ and $\lambda(\min A_m)=-1$,
        \item[(c)] $m>1$ and $\lambda(\min A_m)=1$,
        \item[(d)] $m>1$ and $\lambda(\min A_m)=-1$.
    \end{itemize}
    Notice that if Case (a) applies, the direction of all edges can be reversed and it results in another spanning tree $T_{\Delta^\prime}$ with $\Delta^\prime\in\triang_{(i)}$ for which Case (b) applies and vice versa.
    The same holds for Cases (c) and (d).
    Thus, we get $\inedge(T_{\Delta^\prime})=a-1-\inedge{T_\Delta}$, which means that when we count the elements in $\triang_{(i)}$ with the number of their respective ingoing edges, we can fix without loss of generality the value for $\lambda(\min A_m)$.

    $\mathbf{m=1:}$ We choose $\lambda(\min A_1)=\lambda(r)=1$.
    Let $i$ denote the number of vertices in $P\setminus N$ and let $j$ denote the number of vertices in $A_1^-$.
    The number of ingoing edges in this situation is $i+1+j$.
    This gives rise to the first line of the formula $h^{(i)}$: we sum over all possible choices of $i$ and $j$ and multiply
    the number of ways to pick $P\setminus N$,
    the number of ways to pick $A_1^-$,
    the number of planar spanning trees between $P$ and $A_1^+$,
    the number of planar spanning trees between $N$ and $A_1^-$, and
    polynomial $t^{i+j+1}+t^{a-i-j-2}$, which counts the number of ingoing edges in Cases (a) and (b).
    
    $\mathbf{m>1:}$ We choose $\lambda(\min A_m)=1$.
    Again, we let $i$ denote the number of vertices in $P\setminus N$ and let $j$ denote the number of vertices in $A_1^-$.
    In this case the number of ingoing edges is $i+j$ because $r$ itself is an element of $P$ now and thus cannot be counted in.
    By an analogous statement to the on in the previous case and by summing over all the $A_m$ with $m>1$, we get the second line of the formula which concludes the proof.
\end{proof}

\subsection{Simplices in type (ii) facets}

\begin{figure}
    \centering
    \begin{tikzpicture}
        \foreach \x in {0,1}
            \foreach \y in {0,1,2} 
                \foreach \z in {0,1,2,3,4}
                    \draw[fill] (\x+2.5*\z,\y) circle (0pt) coordinate (n_\x_\y_\z);
        \foreach \x in {0,1}
            \foreach \y in {0,1,2} 
                \foreach \z in {5,6,7,8,9}
                    \draw[fill] (\x+2.5*\z-2.5*5,\y-3.5) circle (0pt) coordinate (n_\x_\y_\z);
        \foreach \x in {0,1}
            \foreach \y in {0,1,2}
                \foreach \z in {10,11,12}
                    \draw[fill] (\x+2.5*\z-2.5*10,\y-7) circle (0pt) coordinate (n_\x_\y_\z);

        \node[shape=circle,draw=black] (c1_a1_1) at (n_0_0_0) {\tiny$I_1$};
        \node[shape=circle,draw=black] (c1_a2_0) at (n_1_1_0) {\tiny$O_2$};
        \node[shape=circle,draw=black] (c1_a3_0) at (n_1_2_0) {\tiny$O_3$};

        \path [->](c1_a2_0) edge node[left] {} (c1_a1_1);
        \path [->](c1_a3_0) edge node[left] {} (c1_a1_1);

        \node[shape=circle,draw=black] (c2_a1_1) at (n_0_0_1) {\tiny$I_1$};
        \node[shape=circle,draw=black] (c2_a2_0) at (n_1_1_1) {\tiny$O_2$};
        \node[shape=circle,draw=black] (c2_a3_1) at (n_0_2_1) {\tiny$I_3$};

        \path [->](c2_a2_0) edge node[left] {} (c2_a1_1);
        \path [->](c2_a2_0) edge node[left] {} (c2_a3_1);

        \node[shape=circle,draw=black] (c3_a1_1) at (n_0_0_2) {\tiny$I_1$};
        \node[shape=circle,draw=black] (c3_a2_1) at (n_0_1_2) {\tiny$I_2$};
        \node[shape=circle,draw=black] (c3_a3_0) at (n_1_2_2) {\tiny$O_3$};

        \path [->](c3_a3_0) edge node[left] {} (c3_a1_1);
        \path [->](c3_a3_0) edge node[left] {} (c3_a2_1);

        \node[shape=circle,draw=black] (c4_a1_1) at (n_0_0_3) {\tiny$I_1$};
        \node[shape=circle,draw=black] (c4_a2_0) at (n_1_1_3) {\tiny$O_2$};
        \node[shape=circle,draw=black] (c4_a3_1) at (n_0_2_3) {\tiny$I_3$};
        \node[shape=circle,draw=black] (c4_a3_0) at (n_1_2_3) {\tiny$O_3$};

        \path [->](c4_a2_0) edge node[left] {} (c4_a1_1);
        \path [->](c4_a3_0) edge node[left] {} (c4_a1_1);
        \path [->](c4_a2_0) edge node[left] {} (c4_a3_1);

        \node[shape=circle,draw=black] (c5_a1_1) at (n_0_0_4) {\tiny$I_1$};
        \node[shape=circle,draw=black] (c5_a2_1) at (n_0_1_4) {\tiny$I_2$};
        \node[shape=circle,draw=black] (c5_a2_0) at (n_1_1_4) {\tiny$O_2$};
        \node[shape=circle,draw=black] (c5_a3_0) at (n_1_2_4) {\tiny$O_3$};

        \path [->](c5_a2_0) edge node[left] {} (c5_a1_1);
        \path [->](c5_a3_0) edge node[left] {} (c5_a1_1);
        \path [->](c5_a3_0) edge node[left] {} (c5_a2_1);

        \node[shape=circle,draw=black] (c6_a1_1) at (n_0_0_5) {\tiny$I_1$};
        \node[shape=circle,draw=black] (c6_a1_0) at (n_1_0_5) {\tiny$O_1$};
        \node[shape=circle,draw=black] (c6_a2_0) at (n_1_1_5) {\tiny$O_2$};
        \node[shape=circle,draw=black] (c6_a3_1) at (n_0_2_5) {\tiny$I_3$};

        \path [->](c6_a2_0) edge node[left] {} (c6_a1_1);
        \path [->](c6_a2_0) edge node[left] {} (c6_a3_1);
        \path [->](c6_a1_0) edge node[left] {} (c6_a3_1);

        \node[shape=circle,draw=black] (c7_a1_1) at (n_0_0_6) {\tiny$I_1$};
        \node[shape=circle,draw=black] (c7_a1_0) at (n_1_0_6) {\tiny$O_1$};
        \node[shape=circle,draw=black] (c7_a2_1) at (n_0_1_6) {\tiny$I_2$};
        \node[shape=circle,draw=black] (c7_a3_0) at (n_1_2_6) {\tiny$O_3$};

        \path [->](c7_a3_0) edge node[left] {} (c7_a1_1);
        \path [->](c7_a3_0) edge node[left] {} (c7_a2_1);
        \path [->](c7_a1_0) edge node[left] {} (c7_a2_1);

        \node[shape=circle,draw=black] (c8_a1_1) at (n_0_0_7) {\tiny$I_1$};
        \node[shape=circle,draw=black] (c8_a2_1) at (n_0_1_7) {\tiny$I_2$};
        \node[shape=circle,draw=black] (c8_a2_0) at (n_1_1_7) {\tiny$O_2$};
        \node[shape=circle,draw=black] (c8_a3_1) at (n_0_2_7) {\tiny$I_3$};
        \node[shape=circle,draw=black] (c8_a3_0) at (n_1_2_7) {\tiny$O_3$};

        \path [->](c8_a2_0) edge node[left] {} (c8_a1_1);
        \path [->](c8_a2_0) edge node[left] {} (c8_a3_1);
        \path [->](c8_a3_0) edge node[left] {} (c8_a1_1);
        \path [->](c8_a3_0) edge node[left] {} (c8_a2_1);

        \node[shape=circle,draw=black] (c9_a1_1) at (n_0_0_8) {\tiny$I_1$};
        \node[shape=circle,draw=black] (c9_a1_0) at (n_1_0_8) {\tiny$O_1$};
        \node[shape=circle,draw=black] (c9_a2_1) at (n_0_1_8) {\tiny$I_2$};
        \node[shape=circle,draw=black] (c9_a3_1) at (n_0_2_8) {\tiny$I_3$};
        \node[shape=circle,draw=black] (c9_a3_0) at (n_1_2_8) {\tiny$O_3$};

        \path [->](c9_a1_0) edge node[left] {} (c9_a2_1);
        \path [->](c9_a1_0) edge node[left] {} (c9_a3_1);
        \path [->](c9_a3_0) edge node[left] {} (c9_a1_1);
        \path [->](c9_a3_0) edge node[left] {} (c9_a2_1);

        \node[shape=circle,draw=black] (c10_a1_1) at (n_0_0_9) {\tiny$I_1$};
        \node[shape=circle,draw=black] (c10_a1_0) at (n_1_0_9) {\tiny$O_1$};
        \node[shape=circle,draw=black] (c10_a2_0) at (n_1_1_9) {\tiny$O_2$};
        \node[shape=circle,draw=black] (c10_a3_1) at (n_0_2_9) {\tiny$I_3$};
        \node[shape=circle,draw=black] (c10_a3_0) at (n_1_2_9) {\tiny$O_3$};

        \path [->](c10_a1_0) edge node[left] {} (c10_a3_1);
        \path [->](c10_a2_0) edge node[left] {} (c10_a1_1);
        \path [->](c10_a2_0) edge node[left] {} (c10_a3_1);
        \path [->](c10_a3_0) edge node[left] {} (c10_a1_1);

        \node[shape=circle,draw=black] (c11_a1_1) at (n_0_0_10) {\tiny$I_1$};
        \node[shape=circle,draw=black] (c11_a1_0) at (n_1_0_10) {\tiny$O_1$};
        \node[shape=circle,draw=black] (c11_a2_1) at (n_0_1_10) {\tiny$I_2$};
        \node[shape=circle,draw=black] (c11_a2_0) at (n_1_1_10) {\tiny$O_2$};
        \node[shape=circle,draw=black] (c11_a3_1) at (n_0_2_10) {\tiny$I_3$};

        \path [->](c11_a1_0) edge node[left] {} (c11_a2_1);
        \path [->](c11_a1_0) edge node[left] {} (c11_a3_1);
        \path [->](c11_a2_0) edge node[left] {} (c11_a1_1);
        \path [->](c11_a2_0) edge node[left] {} (c11_a3_1);

        \node[shape=circle,draw=black] (c12_a1_1) at (n_0_0_11) {\tiny$I_1$};
        \node[shape=circle,draw=black] (c12_a1_0) at (n_1_0_11) {\tiny$O_1$};
        \node[shape=circle,draw=black] (c12_a2_1) at (n_0_1_11) {\tiny$I_2$};
        \node[shape=circle,draw=black] (c12_a2_0) at (n_1_1_11) {\tiny$O_2$};
        \node[shape=circle,draw=black] (c12_a3_0) at (n_1_2_11) {\tiny$O_3$};

        \path [->](c12_a1_0) edge node[left] {} (c12_a2_1);
        \path [->](c12_a2_0) edge node[left] {} (c12_a1_1);
        \path [->](c12_a3_0) edge node[left] {} (c12_a1_1);
        \path [->](c12_a3_0) edge node[left] {} (c12_a2_1);

        \node[shape=circle,draw=black] (c13_a1_1) at (n_0_0_12) {\tiny$I_1$};
        \node[shape=circle,draw=black] (c13_a1_0) at (n_1_0_12) {\tiny$O_1$};
        \node[shape=circle,draw=black] (c13_a2_1) at (n_0_1_12) {\tiny$I_2$};
        \node[shape=circle,draw=black] (c13_a2_0) at (n_1_1_12) {\tiny$O_2$};
        \node[shape=circle,draw=black] (c13_a3_1) at (n_0_2_12) {\tiny$I_3$};
        \node[shape=circle,draw=black] (c13_a3_0) at (n_1_2_12) {\tiny$O_3$};

        \path [->](c13_a1_0) edge node[left] {} (c13_a2_1);
        \path [->](c13_a1_0) edge node[left] {} (c13_a3_1);
        \path [->](c13_a2_0) edge node[left] {} (c13_a1_1);
        \path [->](c13_a2_0) edge node[left] {} (c13_a3_1);
        \path [->](c13_a3_0) edge node[left] {} (c13_a1_1);
        \path [->](c13_a3_0) edge node[left] {} (c13_a2_1);
    \end{tikzpicture}
    \caption{The 13 types of facet graphs of $K_{a,b,c}$.}\label{figure 1}
\end{figure}
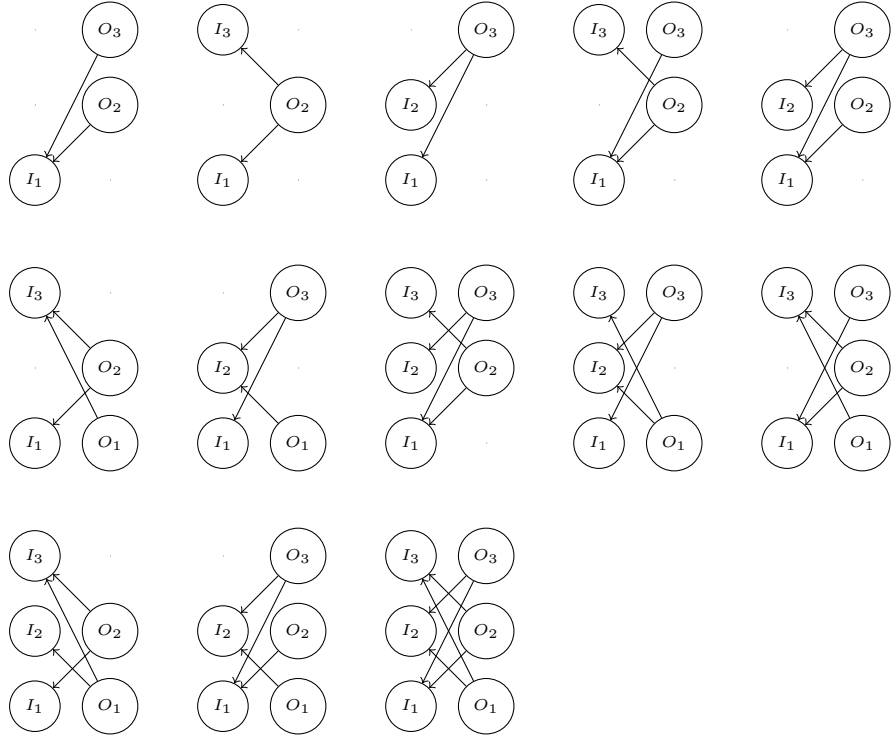

The situation for type (ii) facets is more complicated.
To make things easier, we first define a \emph{labeling in normal form} to be any facet-inducing labeling $\nu$ such that for every multipartite class $A_i$ of $K_{a_1,a_2,\ldots,a_k}$ and every $v,w\in A_i$ that $\nu(v)=1$ and $\nu(w)=0$ implies $v<w$.
Further, define the \emph{opposite labeling} of a facet-inducing $\lambda$ to be the labeling $\Bar{\lambda}$ such that the edge set of $\restr{K_{a_1,a_2,\ldots,a_k}}{\Bar{\lambda}}$ consists of all the reversed edges of $\restr{K_{a_1,a_2,\ldots,a_k}}{\lambda}$.
In the case of type (ii) facets, that means that $\Bar{\lambda}(v) = 1-\lambda(v)$ vor every $v$.
Lastly, for any facet-inducing labeling $\lambda$, the \emph{associated labeling in nomral form} is the labeling in normal form $\nu_\lambda$ such that for every multipartite class of vertices $A_i$, $\abs{\restr{\lambda}{A_i}^{-1}(1)}=\abs{\restr{\nu_\lambda}{A_i}^{-1}(1)}$.

We collect some facts about these objects.

\begin{lemma}
    Let $\lambda$ be a facet-inducing labeling of a type (ii) facet and let $\Delta,\Delta^\prime\in\triang_\lambda$ be simplices.
    Further, let $r$ denote the smallest vertex in $K_{a_1,a_2,\dots,a_k}$.
    The following statements hold.
    \begin{itemize}
        \item[(a)] $\inedge (T_{\Delta}) = \inedge (T_{\Delta^\prime})$. Because of this, we write $\inedge(\lambda)$ to refer to the number of ingoing edges of the spanning trees of the simplices in $\triang_\lambda$.
        \item[(b)] $\inedge(\Bar{\lambda}) = a-\inedge(\lambda-1)$ where $a=a_1+a_2+\dots+a_n$.
        \item[(c)] $\inedge(\nu) = \abs{\nu^{-1}(0)}$ where $\nu$ denotes a labeling in normal form.
        \item[(d)] The number of simplices in $\mathcal{F}_\lambda$ is equal to the number of simplices in $\mathcal{F}_{\nu_\lambda}$ and $\mathcal{F}_{\Bar{\nu}_\lambda}$.
        \item[(e)] Then $\inedge(\lambda)=\inedge(\nu_{\lambda})$ if $r$ is $1$-labeled and $\inedge\lambda=\inedge(\Bar{\nu}_{\Bar{\lambda}})$ if $r$ is $0$-labeled.
    \end{itemize}
\end{lemma}

\begin{proof}
    The crucial insight is the fact that for any given simplex $\Delta\in\triang_\lambda$, the number of ingoing edges does not depend on $\Delta$: $\inedge(T_\Delta)=\abs{\lambda^{-1}(0))\setminus\set{r}}$.
    This can be easily observed by considering that every edge of a $0$-labeled vertex points away from it -- and if it is the one connecting it to the part of the tree which contains $r$, that is an ingoing while the others are outgoing.
    (a) and (c) are immediate corollaries of this.
    (b) is true because by reversing every edge, the outgoing edges become ingoing and vice versa.
    For (d) we may notice that a permutation $\pi$ of the vertices of $K_{a_1,a_2,\dots,a_k}$ within the multipartite classes $A_i$ induces a map from facets to facets.
    If, in addition, we make sure that for all $\ell$-labeled vertices with $\ell\in\set{0,1}$, $w < v$ implies that $\pi(w) < \pi(v)$, then $\pi$ induces a mapping of the simplices in $\triang_{(ii)}$.
    For (e), we get two cases: $\lambda(r)=1$ and $\lambda(r)=0$.
    By default, if $A_1$ contains a single $1$-labeled vertex, $r$ will be $1$-labeled under $\nu_\lambda$.
    Thus, if $\lambda(r)=1$, the first half of the statement follows from (c).
    If $\lambda(r)=0$, $\Bar{\lambda}(r)=1$ and $\inedge(\Bar{\lambda})=\inedge(\nu_{\Bar{\lambda}})$.
    Thus, $\Bar{\nu}_{\Bar{\lambda}}(r)=0$ and the statement follows.
\end{proof}

Now we define $q(\nu)=\abs{\set{\lambda\colon\nu=\nu_\lambda,\ \lambda(r)=1}}$ and $r(\nu)=\abs{\triang_\nu}$.
Since $\nu$ is uniquely defined by the number of $0$-labeled vertices in every multipartite class, we can identify it with the tuple $(\nu_1,\nu_2,\dots,\nu_k)$ which readily gives us a formula for $q$:
\[
    q(\nu_1,\nu_2,\dots,\nu_k)=\prod_{i=1}^k\binom{a_i-\delta_{1,i}}{\nu_i}
\]
where $\delta_{1,i}$ denotes the Kronecker delta, whose function here is to exclude $r$ for the choice.
With the previous lemma, we get
\[
    h^{(ii)}_{a_1,a_2,\ldots,a_k}=\sum_{\Delta\in\triang_{(ii)}} t^{\inedge(T_\Delta)}
    = \sum_{\text{labelings in}\atop\text{normal form $\nu$}} q(\nu)r(\nu)(t^{\nu_1+\nu_2+\dots+\nu_k}+t^{a-1-\nu_1-\nu_2-\dots-\nu_k})
\]
where $a=a_1+\dots+a_k$ again.

To understand $r$, some more work is necessary.
In particular, we will restrict ourselves to the tripartite case.
Figure \ref{figure 1} shows the $13$ different types of facet graphs that are possible.
We call these graphs \emph{class graphs}.
Its vertices and edges are called \emph{class vertices} and \emph{class edges}.
Every class vertex named $O_j$ (resp. $I_j$) represents the set of $0$-marked (resp. $1$-marked) vertices in the corresponding layer.
Every class edge represents the edges of the directed complete bipartite graph between the two corresponding classes of vertices.

Next, we investigate what class graphs tell us about spanning trees corresponding to unimodular simplices.
Firstly, we notice that not every class edge can contain edges in such a spanning tree.
The Gröbner basis elements of type (5) give configurations with involve $6$ class vertices and $3$ class edges.
More precisely, the class edges $\set{O_1,I_2}$, $\set{O_2,I_3}$, and $\set{O_3,I_1}$ cannot all be non-empty at the same time, which turns the last class graph in Figure \ref{figure 1} into three reduced class graphs, each of them missing one of these class edges.

Let us now assume a reduced class graph.
Let $A$ be a class vertex connected to class vertices $B$ and $C$.
Without loss of generality, assume that for every $b\in B$ and every $c\in C$, $b<c$.
The Gröbner basis elements of type (3b) imply that for two distinct vertices $a,a^\prime\in A$ where $a$ is connected to a vertex in $B$ and $a^\prime$ is connected to a vertex in $C$, $a < a^\prime$.
In particular there can be only one vertex in $\Hat{a}\in A$ which connects to both classes.
We denote the set of vertices in $A$ which connect to $B$ (resp. $C$) by $A_B$ (resp. $A_C$).
Analogously, we define the subsets $B_A$ and $C_A$ of vertices which connect to $A$.
Thus we end up with planar spanning trees between the sets $B_A$ and $A_B$ as well as $C_A$ and $A_C$ respectively.

Notice that all reduced class graphs are paths of length $3$, $4$, or $5$.
Thus, consider a class path with vertices $C_1,C_2,\dots,C_n$ of sizes $c_1,c_2,\dots,c_n$.
We define the following function
\begin{align*}
    &c\left(c_1,c_2,\dots,c_n\right) \\
    &= \sum_{j_2=0}^{c_2-1}\sum_{j_3=0}^{c_3-1}\cdots\sum_{j_{n-1}=0}^{c_{n-1}}\binom{c_1+j_2-1}{j_2}\binom{c_{n-1}-j_{n-1}+c_n-2}{c_n-1}\prod_{i=2}^{n-2}\binom{c_i-j_i+j_{i+1}-1}{j_{i+1}}.
\end{align*}
Although this formula looks complicated, its function is fairly straightforward:
Every binomial coefficient $\binom{c_i-j_i+j_{i+1}-1}{j_{i+1}}$ counts the number of planar spanning trees between the sets ${C_i}_{C_{i+1}}$ and ${C_{i+1}}_{C_i}$ where the former has cardinality $j_i+1$ and the latter has cardinality $c_{i+1}-j_{i+1}$.
With this, we get a formula for $r$ in the case of complete tripartite graphs $K_{a,b,c}$.
\[
    r(\nu_1,\nu_2,\nu_3) = \begin{cases}
        c(b,a,c) & \nu_1=0,\nu_2=b,\nu_3=c \\
        c(a,b,c) & \nu_1=0,\nu_2=b,\nu_3=0 \\
        c(a,c,b) & \nu_1=0,\nu_2=0,\nu_3=c \\
        c(\nu_3,a,b,c-\nu_3) & \nu_1=0,\nu_2=b,\nu_3\not\in\set{0,c} \\
        c(\nu_2,a,c,b-\nu_2) & \nu_1=0,\nu_2\not\in\set{0,b},\nu_3=c \\
        c(\nu_1,c,b,a-\nu_1) & \nu_1\not\in\set{0,a},\nu_2=b,\nu_3=0 \\
        c(\nu_1,b,c,a-\nu_1) & \nu_1\not\in\set{0,a},\nu_2=0,\nu_3=c \\
        c(b-\nu_2,\nu_3,a,\nu_2,c-\nu_3) & \nu_1=0,\nu_2\not\in\set{0,b},\nu_3\not\in\set{0,c} \\
        c(a-\nu_1,\nu_3,b,\nu_1,c-\nu_3) & \nu_1\not\in\set{0,a},\nu_2=0,\nu_3\not\in\set{0,c} \\
        c(\nu_1,c-\nu_3,b,a-\nu_1,\nu_3) & \nu_1\not\in\set{0,a},\nu_2=b,\nu_3\not\in\set{0,c} \\
        c(a-\nu_1,\nu_2,c,\nu_1,b-\nu_2) & \nu_1\not\in\set{0,a},\nu_2\not\in\set{0,b},\nu_3=0 \\
        c(\nu_1,b-\nu_2,c,a-\nu_1,\nu_2) & \nu_1\not\in\set{0,a},\nu_2\not\in\set{0,b},\nu_3=c \\
        c(a-\nu_1,\nu_2,c-\nu_3,\nu_1,b-\nu_2,\nu_3) &\\
        + c(\nu_1,c-\nu_3,\nu_2,a-\nu_1,\nu_3,b-\nu_2) &\\
        + c(b-\nu_2,\nu_1,\nu_3,\nu_2,a-\nu_1,c-\nu_3) & \nu_1\not\in\set{0,a},\nu_2\not\in\set{0,b},\nu_3\not\in\set{0,c} \\
        0 & \text{otherwise} \\
    \end{cases}
\]

We can finally assemble the $h^*$-polynomial of $P_{K_{a,b,c}}$.

\begin{thm}\label{theorem: tripartite h*}
    The $h^*$-polynomial of the complete tripartite graph $K_{a,b,c}$ is given by
    \[h^*_{a,b,c}(t) = h_{a,b,c}^{(i)} + h^{(ii)}_{a,b,c}(t).\]
    Here, $h_{a,b,c}^{(i)}$ is given by
    \begin{align*}
        h^{(i)}_{a,b,c}(t) &= \sum_{i=0}^{b+c-1} \sum_{j=1}^{a-1} p(a+b+c,a,i,j)\binom{a+j-i-2}{j-1}\left(t^{i+j+1}+t^{a+b+c-i-j-2}\right) \\
        &+ \sum_{i=0}^{a+c-1} \sum_{j=1}^{b-1} p(a+b+c,b,i,j)\binom{a+c+i-j-1}{a+c-i-1}\left(t^{i+j}+t^{a+b+c-i-j-1}\right) \\
        &+ \sum_{i=0}^{a+b-1} \sum_{j=1}^{c-1} p(a+b+c,c,i,j)\binom{a+b+i-j-1}{a+b-i-1}\left(t^{i+j}+t^{a+b+c-i-j-1}\right)
    \end{align*}
    with
    \[p(x,y,i,j)=\binom{x-y-1}{i}\binom{y-1}{j}\binom{y+i-j-1}{i},\]
    and $h_{a,b,c}^{(ii)}$ is given by
    \begin{align*}
        h^{(ii)}_{a,b,c}(t) &= \sum_{\nu_1=0}^{a-1}\sum_{\nu_2=0}^b\sum_{\nu_3=0}^c
        q(\nu_1,\nu_2,\nu_3)
        r(\nu_1,\nu_2,\nu_3)
        \left(t^{\nu_1+\nu_2+\nu_3}+t^{a+b+c-1-\nu_1-\nu_2-\nu_3}\right)
    \end{align*}
    with
    \[q(\nu_1,\nu_2,\nu_3)=\binom{a-1}{\nu_1}\binom{b}{\nu_2}\binom{c}{\nu_3}\]
    and $r(\nu_1,\nu_2,\nu_3)$ as stated above.
\end{thm}

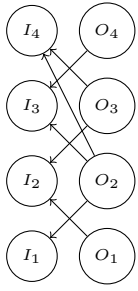
\begin{figure}
    \centering
    \begin{tikzpicture}
        \foreach \x in {0,1}
            \foreach \y in {0,1,2,3}
                \draw[fill] (\x,\y) circle (0pt) coordinate (n_\x_\y);

        \node[shape=circle,draw=black] (a1_1) at (n_0_0) {\tiny$I_1$};
        \node[shape=circle,draw=black] (a1_0) at (n_1_0) {\tiny$O_1$};
        \node[shape=circle,draw=black] (a2_1) at (n_0_1) {\tiny$I_2$};
        \node[shape=circle,draw=black] (a2_0) at (n_1_1) {\tiny$O_2$};
        \node[shape=circle,draw=black] (a3_1) at (n_0_2) {\tiny$I_3$};
        \node[shape=circle,draw=black] (a3_0) at (n_1_2) {\tiny$O_3$};
        \node[shape=circle,draw=black] (a4_1) at (n_0_3) {\tiny$I_4$};
        \node[shape=circle,draw=black] (a4_0) at (n_1_3) {\tiny$O_4$};

        \path [->](a1_0) edge node[left] {} (a2_1);
        \path [->](a2_0) edge node[left] {} (a1_1);
        \path [->](a2_0) edge node[left] {} (a3_1);
        \path [->](a2_0) edge node[left] {} (a4_1);
        \path [->](a3_0) edge node[left] {} (a2_1);
        \path [->](a3_0) edge node[left] {} (a4_1);
        \path [->](a4_0) edge node[left] {} (a3_1);
    \end{tikzpicture}
    \caption{The $O_2$-labeled class vertex has degree $3$.}\label{figure 2}
\end{figure}

For general complete multipartite graphs, the number of possible class graphs grows rapidly as $k$ grows.
Furthermore, in the complete tetrapartite case, class graphs which are not paths start to appear, see e.g. Figure \ref{figure 2}.

\section{New Recursive Relations}

In this section, we gather new evidence for Conjecture \ref{conjecture: HKM}.
First, we state the relevant $h^*$-polynomials.

\begin{prop}\label{proposition: specific h*}
    The $h^*$-polynomials of graphs of the form $K_{1,m,n}$, $K_{1,1,1,n}$, and $K_{2,2,n}$, are given as follows.
    \begin{itemize}
        \item[(a)] $h^*_{1,m,n}(t) = \sum_{i=0}^{\min(m,n)} \binom{2i}{i}\binom{m}{i}\binom{n}{i} t^i(1+t)^{m+n-2i}$
        \item[(b)] $h^*_{1,1,1,n}(t) = 3(n-1)n(1+t)^{n-2}t^2 + 2(2n+1)(1+t)^{n}t + (1+t)^{n+2}$
        \item[(c)] $h^*_{2,2,n}(t) = 20\binom{n}{3}(1+t)^{n-3}t^3 + 2\binom{3n}{2}(1+t)^{n-1}t^2 + 2\binom{3n+1}{1}(1+t)^{n+1}t + (1+t)^{n+3}$
    \end{itemize}
\end{prop}

\begin{proof}
    Notice that (a) is a direct consequence of Propositions \ref{proposition: HJM h*} and \ref{proposition: OT contraction}.
    For (b), see Example \ref{example: h*_1,1,1,n}.
    For (c), one can use the description of $h^*_{a,b,c}$ of Theorem \ref{theorem: tripartite h*} to derive the coefficients of $h^*_{2,2,n}$.
    Then these coefficients can be checked against those of (c).
    Since this is a very tedious process, the reader may consult the corresponding file on 
    \begin{center}
    \url{https://github.com/maxkoelbl/seps_multipartite_graphs/}.
    \end{center}
    It was programmed with SAGEMATH \cite{sage}.
\end{proof}

In the following we will denote the Ehrhart polynomial of $P_{K_{a_1,\dots,a_k}}$ by $E_{a_1,\dots,a_k}$.

\begin{prop}\label{proposition: new relations}
    The following statements hold.
    \begin{itemize}
        \item[(a)] For every $n\geq 2$ there exist $\alpha,\alpha_0\in\RR_{\geq 0}$ which depend on $n$ and satisfy
        \[E_{1,1,n}(x) = \alpha (2x+1) E_{1,n}(x) + \alpha_0 E_{1,n-1}(x).\]
        \item[(b)] For every $n\geq 2$ there exist $\alpha,\alpha_0,\alpha_1\in\RR_{\geq 0}$ which depend on $n$ and satisfy
        \[E_{1,1,n+1}(x) = \alpha (2x+1) E_{1,1,n}(x) + \alpha_0 E_{1,1,n-1}(x) + \alpha_1 E_{1,n}(x).\]
        \item[(c)] For every $n\geq 2$ there exist $\alpha,\alpha_0,\alpha_1\in\RR_{\geq 0}$ which depend on $n$ and satisfy
        \[E_{1,2,n}(x) = \alpha (2x+1) E_{1,1,n}(x) + \alpha_0 E_{1,1,n-1}(x) + \alpha_1 E_{1,n}(x).\]
        \item[(d)] For every $n\geq 2$ there exist $\alpha,\alpha_0,\alpha_1,\alpha_2\in\RR_{\geq 0}$ which depend on $n$ and satisfy
        \[E_{1,2,n+1}(x) = \alpha (2x+1) E_{1,2,n}(x) + \alpha_0 E_{1,2,n-1}(x) + \alpha_1 E_{1,1,n}(x) + \alpha_2 E_{1,n+1}(x).\]
        \item[(e)] For every $n\geq 2$ there exist $\alpha,\alpha_0,\alpha_1\in\RR_{\geq 0}$ which depend on $n$ satisfy
        \[E_{1,1,1,n}(x) = \alpha (2x+1) E_{1,1,n}(x) + \alpha_0 E_{1,1,n-1}(x) + \alpha_1 E_{1,n}(x).\]
        \item[(f)] For every $n\geq 2$ there exist $\alpha,\alpha_0,\alpha_1,\alpha_2\in\RR_{\geq 0}$ which depend on $n$ and satisfy
        \[E_{4,n}(x) = \alpha (2x+1) E_{3,n}(x) + \alpha_0 E_{3,n-1}(x) + \alpha_1 E_{2,n}(x) + \alpha_2 E_{1,n+1}(x).\]
        \item[(g)] For every $n\geq 2$ there exist $\alpha,\alpha_0,\alpha_1,\alpha_2\in\RR_{\geq 0}$ which depend on $n$ and satisfy
        \[E_{3,n+1}(x) = \alpha (2x+1) E_{3,n}(x) + \alpha_0 E_{3,n-1}(x) + \alpha_1 E_{2,n}(x) + \alpha_2 E_{1,n+1}(x).\]
        \item[(h)] For every $n\geq 2$ there exist $\alpha,\alpha_0,\alpha_1,\alpha_2\in\RR_{\geq 0}$ which depend on $n$ and satisfy
        \[E_{2,2,n}(x) = \alpha (2x+1) E_{1,2,n}(x) + \alpha_0 E_{1,2,n-1}(x) + \alpha_1 E_{1,1,n}(x) + \alpha_2 E_{1,n+1}(x).\]
        \item[(i)] For every $n\geq 2$ there exist $\alpha,\alpha_0,\alpha_1,\alpha_2\in\RR_{\geq 0}$ which depend on $n$ and satisfy
        \[E_{1,3,n}(x) = \alpha (2x+1) E_{1,2,n}(x) + \alpha_0 E_{1,2,n-1}(x) + \alpha_1 E_{1,1,n}(x) + \alpha_2 E_{1,n+1}(x).\]
        \item[(j)] For every $n\geq 2$ there exist $\alpha,\alpha_0,\alpha_1,\alpha_2\in\RR_{\geq 0}$ which depend on $n$ satisfy
        \[E_{1,1,1,n+1}(x) = \alpha (2x+1) E_{1,1,1,n}(x) + \alpha_0 E_{1,1,1,n-1}(x) + \alpha_1 E_{1,1,n}(x) + \alpha_2 E_{1,n+1}.\]
    \end{itemize}
\end{prop}

\begin{proof}
    With the formulas in Propositions \ref{proposition: HJM h*} and \ref{proposition: specific h*}, these relations can be obtained algorithmically\footnote{The code for computing explicitly all the coefficients is also available on
    \begin{center}
    \url{https://github.com/maxkoelbl/seps_multipartite_graphs/}.
    \end{center}
    It was also written using SAGEMATH.}.
    We explain the method of proof using relation (a).
    The proof follows that of Proposition 4.5 in \cite{higashitani_kummer_michalek}.
    Since taking the generating function of a polynomial is a linear operation, addition and scalar multiplication translate immediately to Ehrhart series.
    For (a), we need the Ehrhart series of $E_{1,1,n}$, $(2x+1)E_{1,n}$, and $E_{1,n-1}$.
    Notice that multiplying an Ehrhart polynomial by $x$ corresponds to differentiating its Ehrhart series and then multiplying $t$ to it.
    Since the Ehrhart series of $E_{1,n}$ can be written as $\frac{(1+t)^n}{(1-t)^{n+1}}$,
    we get
    \[\frac{{\left(2 \, nt + t + 1\right)} {\left(t + 1\right)}^{n}}{{\left(t^3 -t^2 -t + 1\right)}{\left(1-t\right)}^{n}}\]
    for the Ehrhart series of $(2x+1)E_{1,n}$.
    Next, we form the equation
    \[1 = \frac{\alpha\sum_{k\geq 0} (2k+1)E_{1,n}(k)t^k + \alpha_0 \sum_{k\geq 0}E_{1,n-1}(k)t^k}{\sum_{k\geq 0} E_{1,1,n}(k)t^k}.\]
    Note that there need not be any solutions for $\alpha$ and $\alpha_0$.
    The right-hand side is a rational function of two polynomials where the numerator polynomial involves $\alpha$ and $\alpha_1$.
    Since the right-hand-side is assumed to be equal to one, obtaining a solution is equivalent to a finding asolution of the system of equations
    \[n_i(\alpha, \alpha_0)=d_i\]
    where $n_i$ is the $i$-th degree coefficient of the numerator polynomial and $d_i$ is the $i$-th degree coefficient of the denominator polynomial.
    Since $n_i$ and $d_i$ both depend on $n$, $\alpha$ and $\alpha_0$ do as well.
    We get $\alpha = \frac{n+2}{2(n+1)}$ and $\alpha_0 = \frac{n}{2(n + 1)}$.
\end{proof}

We can state the main result of this section.

\begin{thm}\label{theorem: new relations}
    The following statements hold.
    \begin{itemize}
        \item[(a)] For every $n\geq 1$, $E_{1,n}$ $\CL$-interlaces $E_{1,1,n}$.
        \item[(b)] For every $n\geq 1$, $E_{1,1,n}$ $\CL$-interlaces $E_{1,1,n+1}$.
        \item[(c)] For every $n\geq 1$, $E_{1,1,n}$ $\CL$-interlaces $E_{1,2,n}$.
        \item[(d)] For every $n\geq 1$, $E_{1,1,n}$ $\CL$-interlaces $E_{1,1,1,n}$.
        \item[(e)] For every $n\geq 1$, $E_{3,n}$ $\CL$-interlaces $E_{4,n}$ if $E_{1,n+1}$ $\CL$-interlaces $E_{3,n}$.
        \item[(f)] For every $n\geq 1$, $E_{1,2,n}$ $\CL$-interlaces $E_{1,3,n}$ if $E_{1,n+1}$ $\CL$-interlaces $E_{1,2,n}$.
        \item[(g)] For every $n\geq 1$, $E_{1,2,n}$ $\CL$-interlaces $E_{2,2,n}$ if $E_{1,n+1}$ $\CL$-interlaces $E_{1,2,n}$.
    \end{itemize}
    In particular, for every $n\geq 1$ $E_{x,y,z,n}$ is a $\CL$-polynomial for $x+y+z\leq 3$ and $x,y,z\geq 0$.
\end{thm}

\begin{proof}
    The six labeled statements in this theorem rest entirely on the recursive relations from Proposition \ref{proposition: new relations}, the relation from Proposition \ref{proposition: cross-polynomial relations}, and Proposition \ref{proposition: HKM interlacing}.
    The concluding statement follows from the labeled statements and Proposition \ref{proposition: HKM known interlacings}.
\end{proof}

\section{Recursive Relations and the $\gamma$-vector}

Looking at the recursive relations in Propositions \ref{proposition: HKM known relations} and \ref{proposition: new relations}, we may notice that as the parameters $a_1,\dots,a_{k-1}$ of the multipartite graphs increase, then so does the complexity of the formulas surrounding them.
The results of this section show that this is not a coincidence.
We will show how the existence of a recursion as well as, to some extent, the number of terms it has, are related with the \emph{$\gamma$-vectors} of the $h^*$-polynomials of all the Ehrhart polynomials involved.

\begin{defi}
    Let $h$ be a palindromic polynomial of degree $d$.
    We define the \emph{$\gamma$-vector} as the polynomial $\sum_{i=0}^{\floor*{\frac{d}{2}}}\gamma_i t^i$ such that $h(t) = \sum_{i=0}^{\floor*{\frac{d}{2}}}\gamma_i(1+t)^{d-2i}t^i$.
    We call the degree of the $\gamma$-vector the \emph{$\gamma$-degree} of $h$.
\end{defi}

\begin{lemma}\label{lemma: gammalemma}
    For every integer $d\geq 1$ and every integer $n\geq 0$, the following equation holds.
    \[\sum_{k\geq 0} \left(\sum_{i=0}^{n} (-1)^i\binom{n}{i}\cross_{d+2(n-i)}(k)\right) t^k = \frac{(1+t)^d (4t)^n}{(1-t)^{d+2n+1}}\]
\end{lemma}

\begin{proof}
    There are two key insights.
    The first is the well-known fact that the generating function of $\cross_d$ is $\frac{(1+t)^d}{(1-t)^{d+1}}$.
    The second is that the generating function of $\cross_2(x) - \cross_0(x)$ can be written as $\frac{4t}{(1-t)^3} = \frac{(1+t)^2}{(1-t)^3}-\frac{1}{1-t}$, which can be checked easily.

    The first insight tells us that for real numbers $c_0,c_1,\dots,c_n$, the generating function of $\sum_{i=0}^n c_i\cross_{i}$ can be written as
    \[\frac{1}{1-t}\sum_{i=0}^n c_i\frac{(1+t)^i}{(1-t)^i}.\]
    Using the second insight tells us that
    \[\frac{(1+t)^d(4t)^n}{(1-t)^{d+2n+1}}=\frac{1}{1-t}\frac{(1+t)^d}{(1-t)^d}\left(\frac{(1+t)^2}{(1-t)^2}-1\right)^n.\]
    Finally, with the binomial theorem, we get
    \[\frac{1}{1-t}\sum_{i=0}^n(-1)^i\binom{n}{i}\frac{(1+t)^{d+2(n-i)}}{(1-t)^{d+2(n-i)}},\]
    which concludes the proof.
\end{proof}

\begin{prop}\label{proposition: gamma is cross}
    Let $p$ be a polynomial of degree $d$ and let $h$ be a polynomial defined by
    \[h(t)=(1-t)^{d+1}\sum_{k\geq 0} p(k) t^k.\]
    If $h$ is a palindromic polynomial with $\gamma$-vector $\gamma$, we get
    \[p(x) = \sum_{i=0}^{\deg\gamma}(-1)^i c_i\cross_{d-2i}(x).\]
    where $c_i = \sum_{j=i}^{\deg\gamma}\frac{1}{4^j}\binom{j}{i}\gamma_j$.
    We call the polynomial $\sum_{i=0}^{\deg\gamma}(-1)^i c_i x^i$ the \emph{cross-polynomial coefficients} of $p$.
\end{prop}

\begin{proof}
    We rewrite the generating function of $p$.
    \[
        \frac{h(t)}{(1+t)^{d+1}} = \frac{\sum_{i=0}^d\gamma_i\frac{1}{4^i}(1+t)^{d-2i}(4t)^i}{(1+t)^{d+1}}
    \]
    Splitting up the sum and applying Lemma \ref{lemma: gammalemma}, we get
    \begin{align*}
        \sum_{k\geq 0} \gamma_0\cross_d(k) t^k 
        +\sum_{k\geq 0} \left(\frac{\gamma_0}{4}\cross_d(k)
        -\frac{\gamma_1}{4^2}\cross_{d-2}(k)\right) t^k 
        +\dots +\sum_{k\geq 0} \left(\sum_{i=0}^{n} (-1)^i\frac{\gamma_n}{4^n}\binom{n}{i}\cross_{d+2(n-i)}(k)\right) t^k
    \end{align*}
    Rearranging to sort the sum by the $\cross_i$ yields the claim.
\end{proof}

In the setting of Proposition \ref{proposition: gamma is cross}, we call the $\gamma$-degree of $h$ the \emph{cross-degree} of $p$.

\begin{thm}\label{theorem: recursive relations general}
    Let $f$ be a degree $d+1$ polynomial with cross-degree $m+1$, let $g$ be a degree $d$ polynomial with cross-degree $m$, and let $h_i$ be degree $d-1$ polynomials with cross degree $i$ for $1\leq i\leq m$.
    Then there exist real numbers $\alpha,\alpha_1,\alpha_2,\dots,\alpha_m$ which satisfy
    \[f(x)=(2x+1)\alpha g(x) + \sum_{i=1}^m\alpha_i h_i(x).\]
\end{thm}

\begin{proof}
    Using Proposition \ref{proposition: cross-polynomial relations}, we can see that the degree $d+1$ polynomial $(2x+1) g(x)$ has cross-degree $m+1$. 
    Thus, the right-hand side of the equation can be written as 
    \begin{align*}
        \alpha c_{g,0} & \cross_{d+1} \\
        + (-\alpha c_{(2x+1)g,1}+\alpha_1 c_{h_1,0}+\alpha_2 c_{h_2,0}+\dots+\alpha_m c_{h_m,0}) & \cross_{d-1} \\
        - (-\alpha c_{(2x+1)g,2}+\alpha_2 c_{h_2,1}+\dots+\alpha_m c_{h_m,1})&\cross_{d-3} \\
        &\vdots \\
        +(-1)^m(-\alpha c_{(2x+1)g,m+1}+\alpha_m c_{h_m,m})&\cross_{d-2m+1}
    \end{align*}
    where $c_{(2x+1)g,i}$ is the $i$-th cross-polynomial coefficient of $(2x+1)g(x)$ and the $c_{h_j,i}$ are the cross-polynomial coefficients of the $h_j$.
    This means that in order to get the left-hand side, all we need to do is choose $\alpha,\alpha_m,\alpha_{m-1},\dots,\alpha_1$ in order.
\end{proof}

For complete bipartite graphs, Proposition \ref{proposition: HJM h*} shows that the $\gamma$-degree of the $h^*$-polynomial of $K_{m,n}$ is $\min\set{m,n}-1$.
Thus, we get the following two immediate corollaries.

\begin{cor}\label{corollary: recursion cor 1}
    Let $n$ be a positive integer.
    For $1\leq m\leq n$ there exist $\alpha,\alpha_0,\alpha_1,\dots,\alpha_{m-1}\in\RR$ such that
    \[E_{m+1,n+1}(x)=(2x+1)\alpha E_{m,n+1}(x)+\sum_{i=0}^{m-1}\alpha_i E_{m-i,n+i}(x)\]
    is satisfied, where $E_{m,n}$ denotes the Ehrhart polynomial of the symmetric edge polytope of $K_{m,n}$.
\end{cor}

\begin{cor}\label{corollary: recursion cor 2}
    Let $n$ be a positive integer.
    For $1\leq m\leq n$ there exist $\alpha,\alpha_0,\alpha_1,\dots,\alpha_{m-1}\in\RR$ such that
    \[E_{m,n+1}(x)=(2x+1)\alpha E_{m,n}(x)+\sum_{i=0}^{m-1}\alpha_i E_{m-i,n+i-1}(x)\]
    is satisfied, where $E_{m,n}$ denotes the Ehrhart polynomial of the symmetric edge polytope of $K_{m,n}$.
\end{cor}

\begin{rem}
    These corollaries alone are not enough to prove Conjecture \ref{conjecture: HKM} for all $K_{m,n}$ for two crucial reasons.
    Firstly, as $m$ increases, the number of interlacings having to be satisfied increases too, and they are between polynomials whose cross-degrees puts them outside the scope of Theorem \ref{theorem: recursive relations general}.
    This is noticeable in the last four statements of Theorem \ref{theorem: new relations} where the interlacing of cross-degree $3$ polynomials by cross-degree $2$-polynomials depend on the interlacing of a cross-degree $2$-polynomial by a cross-degree $0$ polynomial.
    
    Secondly, there is no guarantee that the coefficients $\alpha,\alpha_1,\dots,\alpha_m$ are non-negative, although explicit computations for low $m$ in the context of Corollary \ref{corollary: recursion cor 1} always yield positive coefficients.
    In fact, for $m\geq 4$, explicit computations reveal that $\alpha_2,\dots,\alpha_{m-2}$ are always negative.
    In the case $m=4$, we get $\alpha_2=\frac{n - n^3}{8(5n^3 + 39n^2 + 100n + 96)}$.
    To see the parameters for every $1\leq m\leq 10$, we refer once again to the corresponding SAGEMATH code on
    \begin{center}
        \url{https://github.com/maxkoelbl/seps_multipartite_graphs/}.
    \end{center}
\end{rem}

\begin{conj}
    Let $a_1\leq a_2\leq\dots\leq a_k\leq n$ be positive integers and let $m$ denote the cross-degree of the Ehrhart polynomial of the symmetric edge polytope of $K_{a_1,a_2,\dots,a_k}$.
    Then we conjecture the inequalities
    \[\floor*{\frac{\sum_{i=1}^k a_i}{2}}\leq m+1\leq \sum_{i=1}^k a_i.\]
    Furthermore, we conjecture that the Ehrhart polynomial of the symmetric edge polytope of the graph $K_{1^k, n}$ interlaces that of $K_{1^{k+1}, n}$, where $1^k$ represents a list of ones.
\end{conj}

\printbibliography

\end{document}